\documentclass[12pt]{amsart}
\newtheorem{theorem}{Theorem}[section]
\newtheorem{lemma}[theorem]{Lemma}
\newtheorem{corollary}[theorem]{Corollary}
\newtheorem{proposition}[theorem]{Proposition}
\newtheorem{remark}[theorem]{Remark}
\newtheorem{definition}[theorem]{Definition}
\newtheorem{conjecture}[theorem]{Conjecture}
\newtheorem{question}[theorem]{Question}

\newcommand{\ncom}{\newcommand}
\ncom{\rar}{\rightarrow}
\ncom{\lrar}{\longrightarrow}
\ncom{\ov}{\overline}
\ncom{\m}{\mbox}
\ncom{\sta}{\stackrel}
\ncom{\comx}{{\mathbb C}}
\ncom{\Z}{{\mathbb Z}}
\ncom{\Q}{{\mathbb Q}}
\ncom{\R}{{\mathbb R}}
\ncom{\G}{{\mathbb G}}
\ncom{\al}{\alpha}
\ncom{\p}{{\mathbb P}}
\ncom{\E}{{\mathbb E}}
\ncom{\N}{{\mathbb N}}
\ncom{\K}{{\mathbb K}}
\ncom{\Le}{{\mathbb L}}
\ncom{\A}{{\mathbb A}}
\ncom{\F}{{\mathbb F}}

\ncom{\f}{\frac}
\ncom{\cA}{{\mathcal A}}
\ncom{\cX}{{\mathcal X}}
\ncom{\cO}{{\mathcal O}}
\ncom{\cW}{{\mathcal W}}
\ncom{\cL}{{\mathcal L}}
\ncom{\cP}{{\mathcal P}}
\ncom{\cH}{{\mathcal H}}
\ncom{\cS}{{\mathcal S}}
\ncom{\cM}{{\mathcal M}}
\ncom{\cC}{{\mathcal C}}
\ncom{\cT}{{\mathcal T}}
\ncom{\cF}{{\mathcal F}}
\ncom{\cN}{{\mathcal N}}
\ncom{\cJ}{{\mathcal J}}
\ncom{\cV}{{\mathcal V}}
\ncom{\cZ}{{\mathcal Z}}
\ncom{\cU}{{\mathcal U}}
\ncom{\cSU}{{\mathcal S \mathcal U}}
\ncom{\cG}{{\mathcal G}}
\ncom{\cQ}{{\mathcal Q}}
\ncom{\cR}{{\mathcal R}}

\ncom{\eop}{{\hfill $\Box$}}
\begin{document}
\baselineskip=16pt

\title[Murre's conjectures and finite dimensionality]{Murre's conjectures and explicit Chow--K\"unneth projectors for some varieties}

\author[J. N. Iyer]{Jaya NN Iyer}
\address{School of Mathematics, Institute for Advanced Study, 1 Einstein Drive, Princeton, NJ 08540, USA}
\email{jniyer@ias.edu}
\address{The Institute of Mathematical Sciences, CIT
Campus, Taramani, Chennai 600113, India}
\email{jniyer@imsc.res.in}

\footnotetext{Mathematics Classification Number: 14C25, 14D05, 14D20, 14D21 }
\footnotetext{Keywords: Homogenous spaces, Chow groups, projectors.}


\begin{abstract}
In this paper, we investigate Murre's conjectures 
on the structure of rational Chow groups and exhibit explicit Chow--K\"unneth projectors for some examples. More precisely,
the examples we study are the varieties which have a nef tangent bundle. 
For surfaces and threefolds which have a nef tangent bundle explicit Chow--K\"unneth
projectors are obtained which satisfy Murre's conjectures and the motivic Hard Lefschetz theorem is verified.    
\end{abstract}

\maketitle

\section{Introduction}

Suppose $X$ is a nonsingular projective variety of dimension $n$ defined over 
the complex numbers. Let $CH^i(X)\otimes \Q$ be the Chow group of codimension 
$i$ algebraic cyles modulo rational equivalence. There has been a wide interest
in understanding the structure of the Chow 
groups and their study has led to many conjectural filtrations on the Chow groups. Jacob Murre \cite{Mu2}, \cite{Mu3} has made the following conjectures:

(A) The motive $h(X):=(X,\Delta_X)$ of $X$ has a Chow-K\"unneth decomposition:
  $$\Delta_X= \sum_{i=0}^{2n}\pi_i \in CH^n(X\times X)\otimes \Q$$
such that $\pi_i$ are orthogonal projectors (see \S \ref{CK}).
 
(B) The correspondences $\pi_0,\pi_1,...,\pi_{j-1},\pi_{2j+1},...,\pi_{2n}$ act as zero on $CH^j(X)\otimes \Q$.

(C) Suppose $$F^rCH^j(X)\otimes \Q= Ker \pi_{2j}\cap Ker \pi_{2j-1}\cap...\cap Ker \pi_{2j-r+1}.$$
 Then the filtration $F^\bullet$ of $CH^j(X)\otimes \Q$ is independent of the choice of the projectors $\pi_i$.

(D) Further, $F^1CH^i(X)\otimes \Q= (CH^i(X)\otimes \Q)_{hom}$, the cycles which are homologous to zero.

A related conjecture is that the motive of a smooth projective variety satisfies a motivic
Hard-Lefschetz-Theorem (see \S 2.5). 
A part of these conjectures have been verified for curves, surfaces,
a product of a curve and surface \cite{Mu1}, \cite{Mu3}, abelian varieties and abelian schemes \cite{Sh},\cite{De-Mu}, 
uniruled threefolds \cite{dA-Mul},  elliptic modular varieties
\cite{Go-Mu}, \cite{GHM2}), universal families over Picard modular surfaces \cite{MM} and finite group quotients (maybe singular) of 
abelian varieties \cite{Ak-Jo}. A criterion for the existence of a Chow--K\"unneth decomposition is also prescribed in \cite{Sa}. 

On the other hand, S.I. Kimura \cite{Ki} has introduced the notion
of a finite dimensional motive. He showed the finite
dimensionality of the motive of a curve and an abelian variety.
He has also shown that finite dimensionality is preserved under surjective
morphisms, for product varieties and quotients (see \S \ref{Kim}). He has conjectured that any nonsingular projective
variety $X$ has a finite dimensional motive. Furthermore, if $X$
is known to have a K\"unneth decomposition then finite
dimensionality implies the existence of a Chow--K\"unneth
decomposition for $X$ \cite[Corollary 9]{Gu-Pe}. Hence Kimura's conjecture
is related to Murre's conjectures.

Our aim here is to verify Murre's conjectures and exhibit explicit Chow-K\"unneth projectors, in the following situation. This also gives the existence of the conjectural filtration on the rational Chow groups for these varieties.

\begin{theorem}\label{main}
Suppose $X$ is a nonsingular projective surface or a threefold which has a numerically effective tangent bundle. 
Then the motive  of
$X$ is finite dimensional and we have an explicit Chow-K\"unneth decomposition in terms of the projectors of its albanese reduction.
Moreover, Murre's conjectures is fulfilled by these projectors  and the motivic Hard Lefschetz theorem holds .
\end{theorem}
See Corollary \ref{Murre} and  Proposition \ref{threefolds} in \S \ref{projectors} for the explicit projectors.

Examples of varieties which have a nef tangent bundle include abelian varieties, hyperelliptic varieties
and bundles of homogenous varieties over an abelian variiety.

A similar statement as in Theorem \ref{main} holds for higher dimensional ($\geq 4$) nonsingular varieties with a nef tangent bundle if we assume that
 homogenous bundles $Z/S$ have a relative Chow--K\"unneth decomposition (in the sense of \cite{De-Mu}) and any nonsingular Fano variety with a nef  tangent bundle is a homogenous variety, see \S \ref{remarks}. 

We note that our results imply that whenever a nonsingular projective variety $X$ admits a finite surjective 
cover $Z\lrar X$ such that $Z\lrar A$ is a relative cellular variety over an abelian variety $A$, then
$X$ has an absolute Chow--K\"unneth decomposition and the consequences of Theorem \ref{main} hold, see Corollary \ref{cellular}.  

The proof of Theorem \ref{main} involves  the classification results (\cite{Ca-Pe}, \cite{DPS}) and
applying the criterion in \cite{GHM2}, to obtain absolute Chow--K\"unneth projectors from relative Chow--K\"unneth projectors.
 We also write down the explicit Chow--K\"unneth projectors for $X$ in terms of the projectors
of its Albanese reduction (see Proposition \ref{pr.-ckx}). These projectors moreover satisfy Poincar\'e duality and
part of Murre'e conjectures (similar to the results on abelian varieties).
Thus we apply the geometric methods and results on classification theory to answer some questions on algebraic cycles.

Here is an outline of the paper: we briefly review the preliminaries in \S 2. 
A study of the above conjectures for rational homogeneous bundles over 
varieties in some generality followed by 
a more detailed discussion for varieties with nef tangent bundle is carried out in \S 3. 

{\Small Acknowledgements: We are grateful to R. Joshua for asking the question to 
extend some results in \cite{Ak-Jo} to other examples and for the hospitality
 during a visit in Oct. 2005. During the preparation of this article, he 
brought to our attention the motivic Hard Lefschetz theorem and issues on the 
existence of relative cellular decomposition and we are thankful to him. 
We also thank S-I. Kimura for explaining his work \cite{Ki}, S. M\"uller-Stach for his comments on the note and P. Deligne for pointing out some errors in the earlier version.}

\section{Preliminaries}
We work over the field of complex numbers in this paper and all the Chow groups are taken with $\Q$ -coefficients.

\subsection{Category of motives}
The category of nonsingular projective varieties over $\comx$ will be denoted by $\cV$.
For an object $X$ of $\cV$, let $CH^i(X)_\Q=CH^i(X)\otimes \Q$ denote the rational Chow group of codimension $i$ algebraic cycles modulo rational equivalence.
 We will use the standard framework of the category of Chow motives $\cM_{rat}$ in this paper and refer to \cite{Mu2} for details.
We denote the category of motives $\cM_{\sim}$, where $\sim$ is any equivalence, for instance $\sim$ is homological or numerical equivalence.
When $S$ is a smooth variety, we also consider the category of relative 
Chow motives $CH\cM(S)$ which is introduced in \cite{De-Mu} and \cite{GHM}. When $S=\m{Spec } \comx$ then the category $CH\cM(S)= \cM_{rat}$. 

\subsection{Chow--K\"unneth decomposition for a variety}\label{CK}

Suppose $X$ is a nonsingular projective variety over $\comx$ of dimension $n$.
Let $\Delta_X\subset X\times X$ be the diagonal.
Consider the K\"unneth decomposition of $\Delta$ in the Betti Cohomology:
$$\Delta_X = \oplus_{i=0}^{2n}\pi_i^{hom}$$
where $\pi_i^{hom}\in H^{2n-i}(X)\otimes H^i(X)$.

\begin{definition}
The motive of $X$ is said to have K\"unneth decomposition if each of the
classes $\pi_i^{hom}$ are algebraic i.e., $\pi_i^{hom}$ is the image of an algebraic cycle $\pi_i$
under the cycle class map from the rational Chow groups to the Betti Cohomology.

\end{definition}

\begin{definition}
The motive of $X$ is said to have a Chow--K\"unneth decomposition if each of the
classes $\pi_i^{hom}$ is algebraic and are orthogonal projectors, i.e.,
$\pi_i\circ \pi_j=\delta_{i,j}\pi_i$.
\end{definition}

We recall the following well-known statement which we will need in Proposition \ref{pr.-ckx}.
\begin{lemma}\label{le.-2}
Suppose $f:X\lrar Y$ is a finite surjective morphism. If $X$ has a
K\"unneth decomposition then $Y$ also has a K\"unneth decomposition.
\end{lemma}
\begin{proof}
Suppose $\pi_i^X$ are the K\"unneth components for $X$. Since
$f_*\Delta_X=m.\Delta_Y$, where $m=degree (f)$, it follows that
$(1/m).f_*(\pi_i^X)$ are the K\"unneth components for $Y$.
\end{proof}


\subsection{Relative Chow--K\"unneth decomposition}

Suppose $f:X\lrar S$ is a smooth projective morphism of relative dimension $d$.
A relative Chow--K\"unneth decomposition of $X\lrar S$ is a decomposition of the relative diagonal
cycle 
$$
\Delta_S =\sum_{i=0}^{2g}\Pi_i\,\in\,CH^d(X\times_S X)\otimes \Q
$$
such that $\Pi_i\circ \Pi_j= \delta_{i,j}\Pi_i$, i.e., they are orthogonal projectors.

It is conjectured that any smooth projective morphism as above admits a relative Chow--K\"unneth decomposition and instead of the cohomology, the projectors act on $\textbf{R}f_*\Q$ in the derived category of bounded complexes of cohomlogically constructible $\Q$- sheaves on $S$ (see \cite{De-Mu}). 
 More generally, it is  conjectured by Corti and Hanamura \cite{Co-Ha} that a \textit{motivic Decomposition theorem} is true for a stratified projective morphism. 

We recall the criterion in \cite{GHM2} to obtain absolute Chow-K\"unneth projectors from the relative Chow--K\"unneth projectors. We will need this in \S \ref{absmotive}. Since we will be dealing with a smooth projective morphism $\pi: X\rar S$, we state the criterion in \cite{GHM2}
in this case.

We first state the condition to have an explicit decomposition;

\begin{lemma}\label{explicit}
Let $X$ and $S$ be smooth quasi-projective varieties over $\comx$ and $\pi:X\rar S$ is a projective morphism. Then 

a) giving a projector $\Pi\in \m{CH}_{\m{dim }X}(X\times _S X)$ and an isomorphism in $CH\cM(S)$
$$
f: (X,\Pi,0)\lrar \oplus^m(S,\Delta(S),-q)
$$
is equivalent to giving elements $f_1,...,f_m \in \m{CH}_{\m{dim }X+q}(X)$ and $g_1,...,g_m \in \m{CH}_{\m{dim }X-q}(X)$ subject to the condition $\pi_*(f_i.g_j)=\delta_{ij}[S]$. Here $f_i.g_j$ is the intersection product in $X$ and $[S]\in \m{CH}_{\m{dim }S}(S)$ is the fundamental class.

b) if $(X,P,0) \in CHM(S)$ and $\Pi$ is a constituent of $P$ (i.e., $P\circ \Pi=\Pi\circ P=\Pi$) then we can take the above $f_i$ and $g_j$ such that moreover $f_i\circ P=f_i$ and $g_j=P\circ g_j$. Conversely if we have such $f_i$ and $g_j$ then the corresponding $\Pi$ is a constituent of $P$.
\end{lemma}
\begin{proof}
This is \cite[Lemma 4]{GHM2}.
\end{proof}

\begin{theorem}\label{GHM2}
Suppose the variety $\pi:X\rar S$ has  a relative Chow--K\"unneth decomposition. This means that there is an orthogonal set of projectors $P^i\in \m{CH}_{{dim }X}(X\times_S X)$ with $\sum_iP^i=1$ such that $P^i$ acts on $R^j\pi_*\Q_X$ as identity if $j=i$ and as zero if $j\neq i$.
Assume the following conditions;

1. $S$ has a Chow--K\"unneth decomposition over $\comx$.

2. If $t$ is a point of $S$, the natural map
$$
CH^r(X)\lrar H^{2r}_B(X_t(\comx),\Q)^{\pi_1^{top}(S,t)}
$$
is surjective for $0\leq r\leq d:=\m{dim} X/S$. The target is the invariant part under the action of the topological fundamental group of $S$.

3. For $i$ odd, $H^i_B(X_t,\Q)^{\pi_1^{top}(S,t)}=0$.

4. If $\cV$ is the local system $R^i\pi_*\Q_X$, for any $i=2r-1,\,0\leq r\leq d$, or $R^{2r}\pi_*\Q_X/(R^{2r}\pi_*\Q_X)^{\pi_1}$, if $0\leq r\leq d$, then $H^q(S,\cV)=0 $ if $q\neq \m{dim } S$.

Under these assumptions $X$ has an absolute Chow--K\"unneth decomposition over $\comx$.
\end{theorem}
\begin{proof}
This is \cite[Main Theorem 1.3]{GHM2}. In loc.cit the local system in assumption 4 is arbitrary but we notice that in the proof only those local systems as stated above, are required to satisy assumption 4.
\end{proof}

\begin{lemma}\label{alg}
Construction of a projector $(P^{2r}/S)_{alg}$ (the subscript ``alg'' indicates it splits the algebraic
 cohomology in the fibre $X_t$) which is a constituent of a projector $P^{2r}/S$, together with an isomorphism
$$
(X/S,(P^{2r}/S)_{alg},0)\simeq \oplus^m(S,\Delta(S), -r)
$$
follows from assumption 2. of Theorem \ref{GHM2}.
\end{lemma} 
\begin{proof}
This is proved in Step II in the proof of \cite[Main theorem]{GHM2}.
\end{proof}

We will apply these statements when we look at rational homogenous bundles over abelian varieties in \S \ref{absmotive}.

\subsection{Finite dimensional motives}\label{Kim}

The notion of a \textit{finite dimensional motive} was  introduced by
S.I. Kimura . We refer to his paper \cite{Ki} for the definition and properties.
Examples of finite dimensional motives includes the case of curves and abelian varieties.
We recall the following results from \cite{Ki} which we will use in this paper. These provide further examples of finite dimensional motives.

\begin{lemma}\cite[Corollary 5.11]{Ki}.
\label{le.-3}
If the motive of $X$ and $Y$ are finite dimensional then the motive of
the product $X\times Y$, the sum $X\oplus Y$ and a direct summand of a finite dimensional motive, are also finite dimensional.
\end{lemma}

\begin{lemma}\cite[Proposition 6.9]{Ki}.
\label{le.-4}
Suppose $f:X\lrar Y$ is a surjective morphism. If the motive of $X$
is finite dimensional then the motive of $Y$ is also finite dimensional.
\end{lemma}

Kimura (\cite[Proposition 7.5]{Ki}) observed that the finite dimensionality of a motive $M$ implies the nilpotence of the ideal $I$ defined by the exact sequence
$$
0\lrar I \lrar \m{End}_{\cM_{rat}}(M) \lrar \m{End}_{\cM_{hom}}(M)\lrar 0.
$$

A result due to U. Jannsen (\cite[5.3]{Ja}) is used to deduce the following:

\begin{proposition}\cite[Corollary 9]{Gu-Pe}.
\label{pr.-1}
Suppose $h(X)$ is finite dimensional. Assume that the K\"unneth components of
$X$ are algebraic. Then $X$ has a Chow--K\"unneth decomposition.
\end{proposition}

\subsection{Lefschetz decomposition of a motive}\label{lef}
We recall the definitions of Lefschetz operator and Lefschetz decomposition of motives. We refer to \cite{Kl}, \cite{Ku}, for the details.
 
Given a smooth projective variety $X$ of dimension $n$ and an ample divisor $\al\in CH^1(X)_Q$, we can define the Lefschetz operator $L_\al$.
Consider the closed embedding $i: \al\hookrightarrow X$.
Define $$L_\al:= i_*i^*:h(X)\lrar h(X)(1)$$ and the $(n-1)$-fold compositions with itself is
 denoted by
$$L^{n-i}_\al: h(X)\lrar h(X)(n-i).$$

In the realizations $L_\al$ induces the multiplication with
the class of $\al$. 
Then the Hard Lefschetz theorem gives an isomorphism between the cohomology 
groups defined by the Lefschetz operator :
$$L_\al^{n-i}: H^i(X,\Q)\sta{\simeq}{\lrar} H^{2n-i}(X,\Q).$$
Furthermore, there is a decomposition of the cohomology into primitive cohomologies,
$$H^i(X,\Q)=\bigoplus_{\m{max}\{0,i-n\}\leq k\leq [i/2]}L^kP^{i-2k}(X,\Q).$$

The standard conjectures of Lefschetz type predict a motivic statement
of above statements, namely, there is an isomorphism of motives
$$
L_\al^{n-i}: h^i(X)\sta{\simeq}{\lrar} h^{2n-i}(X)(n-i)
$$
and a decomposition of motives 
$$h^i(X)=  \bigoplus_{\m{max}\{0,i-n\}\leq k\leq [i/2]}L_\al^kP^{i-2k}(X).$$
Here $P^i(X)$ is the motivic analogue of the primitive cohomology.

\begin{theorem}\label{th.-KuSc}
Suppose $X$ is an abelian variety then the motivic Hard Lefschetz theorem and the Lefschetz decomposition hold for $X$.
\end{theorem}
\begin{proof}
This was proved by K\"unnemann and Scholl, see \cite{Ku}, \cite{Sc}.
\end{proof} 


\section{Rational homogenous bundles over a variety and varieties with a nef tangent bundle}

In this section we begin by recalling the motive a homogenous space and of a relative cellular variety.
These results will be used to study the motive of varieties with numerically effective tangent bundle.

\subsection{The motive of a rational homogenous space}

Let $G$ be a
 reductive linear algebraic group and $P$ be a
parabolic subgroup of $G$. Then $F:=G/P$ is a complete variety. 
Notice that $F$ is a cellular variety, i.e., it has a cellular decomposition
$$\emptyset=F_{-1}\subset F_0 \subset...\subset F_n=F$$
such that each $F_i\subset F$ is a closed subvariety and $F_i-F_{i-1}$ is an affine space.

Then we have
\begin{lemma}\cite[Theorem, p.363]{Ko}.
\label{le.-motF}
The Chow motive $h(F)=(F,\Delta_F)$ of $F$ decomposes as a direct sum of twisted Tate motives
$$h(F)=\bigoplus_{\omega}\Le^{\otimes \m{dim }\omega}.$$
Here $\omega$ runs over the set of cells of $F$.
\end{lemma}


\subsection{The motive of a rational homogenous bundle}\label{absmotive}

Consider a rational homogenous bundle
$$
f:Z\lrar Y
$$
i.e., $f$ is a smooth projective morphism and any fibre $f^{-1}y$ is a rational homogenous variety $G/P$, for some reductive linear algebraic group $G$ and a
parabolic subgroup $P$ of $G$.
 
If $Z\lrar Y$ is a relative cellular variety then it might be possible to prove the following lemma directly, as mentioned to us by R. Joshua. We consider a more general situation when we have a 
homogenous bundle and under the assumption of existence a relative Chow--K\"unneth decomposition we apply the criterion in \cite{GHM2} to deduce absolute
Chow--K\"unneth decomposition.

\begin{lemma}\label{le.-kd}
Let $Y$ be a smooth variety which has a Chow--K\"unneth decomposition.
Suppose $Z\lrar Y$ is a rational homogenous bundle and $d:=\m{dim }Z-\m{dim }Y$. Assume that the bundle $Z\lrar Y$ has a relative Chow--K\"unneth decomposition in the category of relative Chow motives $CH\cM(Y)$ over $Y$.
Then the absolute motive of $Z$ also admits a Chow--K\"unneth decomposition.
\end{lemma}
\begin{proof}
By assumption, we have relative Chow--K\"unneth projectors for the variety $Z\lrar Y$. Write a relative Chow--K\"unneth decomposition of the relative 
Chow motive of $Z/Y$
$$
(Z/Y,1,0)= \bigoplus_{r=0}^{2d}(Z/Y,\Pi_r,0)
$$
in $CH\cM(Y)$.
We want to apply Theorem \ref{GHM2}, to deduce the absolute Chow--K\"unneth decomposition
from the relative Chow--K\"unneth decomposition.
For this purpose we consider  the ordinary rational cohomology of the  fibres of $f$ and show that the 
cohomology classes of the fibres are invariant under the topological fundamental group of the base $Y$
under the natural monodromy action. More precisely,  the fibres of $f$  are smooth homogenous varieties and have only algebraic cohomology. Hence for any $t\in Y$, it follows that the map
$$
CH^r(Z)\lrar H^{2r}_B(Z_t,\Q)^{\pi^{top}_1(Y,t)}
$$
is surjective, for $0\leq r\leq d$
and for $i$ odd
$$
H^i_B(Z_t,\Q)^{\pi^{top}_1(Y,t)}=0.
$$
Hence assumptions 1.,2., and 3., in Thereom \ref{GHM2} are fulfilled. Moreover the local 
systems $R^i\pi_*\Q_Z$, for any $i=2r-1,\,0\leq r\leq d$, or $R^{2r}\pi_*\Q_Z/(R^{2r}\pi_*\Q_Z)^{\pi_1}$, if $0\leq r\leq d$, are vanishing. The assumption 4. in Thereom
\ref{GHM2} is now also satisfied. Applying Theorem \ref{GHM2}, we obtain an absolute Chow--K\"unneth decomposition for $Z$.

\end{proof}

\begin{corollary}\label{co.-motZ} 
With notations and assumptions as in Lemma \ref{le.-kd},
the motive of $Z$ is
$$
(Z,\Delta_Z)= \bigoplus_i(Z,\pi_i^Z)
$$
where $(Z,\pi^Z_i)= \bigoplus_{j+k=i}m_j.\Le^j\otimes (Y,\pi_k^Y)$.
Here $m_j$ is the number of $j$-codimensional cells on a fibre $\F$ and $(Y,\Delta_Y)=\oplus_k
(Y,\pi^Y_k)$ is a Chow--K\"unneth decomposition for $Y$.
\end{corollary}
\begin{proof}
Now applying Lemma \ref{alg} we get an isomorphism of the 
relative Chow motives 
$$
(Z/Y,\Pi_{2r},0) \simeq (Y,1,-r)^{\oplus_{m_r}}
$$
in $CH\cM(Y)$, where $m_r=\m{dim}H^{2r}_B(Z_t,\Q)$. Using \cite[Lemma 1]{GHM2}
the natural map
$$
\phi:CH_{\m{dim }Z}(Z\times_{Y} Z)\lrar CH_{\m{dim }Z}(Z\times_\comx Z)
$$
is a ring homomorphism and hence transforms orthogonal projectors over $Y$ to
orthogonal projectors over $\comx$. Hence we have a functor
\begin{eqnarray*}
\Phi:CH\cM(Y)& \lrar & \cM_{rat}(\comx) \\
 (Z/Y,\Pi_{2r},i)& \mapsto & (Z,P_{2r},i)
\end{eqnarray*}
Since the cohomology of the fibres $Z_t$ are algebraic, by \cite[\S1.7, \textbf{Step V}]{GHM2}, there is an isomorphism
$$
(Z,P_{2r},0)\simeq (Y,1,-r)^{\oplus m_r}.
$$
Since $Y$ has a Chow-K\"unneth decomposition we can write this decomposition as
$$
(Y,\Delta_Y)= \oplus_l(Y,\pi^Y_l).
$$
This gives
\begin{eqnarray*}
(Z,\Delta_Z)& = & \oplus_{r=0}^d(Z,P_{2r}) \\
            &=& \oplus_{r=0}^d\oplus_k m_r.\Le^r\otimes (Y,\pi_k^Y).
\end{eqnarray*}

Write 
$$
(Z,\pi^Z_i):= \bigoplus_{j+k=i}m_j.\Le^j\otimes (Y,\pi_k^Y)
$$
to give a Chow--K\"unneth decomposition of $(Z,\Delta_Z)= \oplus_i(Z,\pi^Z_i)$.
\end{proof}

\begin{corollary}\label{co.-fd}
With notations and assumptions as in Lemma \ref{le.-kd},
the motive of $Z$ is finite dimensional if the motive of $Y$ is
finite dimensional.
\end{corollary}
\begin{proof}
This follows by Lemma \ref{le.-3} and Corollary \ref{co.-motZ}.
\end{proof}

\begin{corollary}
Suppose $Z\lrar Y$ is a homogenous bundle which has a relative cellular decomposition.
If  the motive of $Y$ is finite dimensional then the motive of $Z$ is also finite dimensional.
\end{corollary}
\begin{proof}
By \cite{Ko} and \cite{Ne-Za}, we know that the motive of $Y$ has a relative Chow-K\"unneth decomposition over $Z$. So we can apply Corollary \ref{co.-motZ} and Corollary \ref{co.-fd}. 
\end{proof}

Fix an ample divisor $\al\in CH^1(Z)_\Q$ which is the pullback of an ample divisor $\beta \in CH^1(Y)_\Q$ twisted by a relatively ample divisor $\gamma \in CH^1(Z)$. Consider the Lefschetz operator $L_\al \in CH^*(Z\times Z)_\Q$, as defined in \S \ref{lef}. Let $m=\m{dim }Y$, $f=\m{ dim }\F$ and $n=m+f$. Here $\F$ is a fibre of $Z\rar Y$.
 
\begin{corollary}\label{co.-Hlef}
Suppose the motivic Hard Lefschetz theorem and Lefschetz decomposition holds
for $Y$. Then it also holds for the variety $Z$.
\end{corollary} 
\begin{proof}
By Corollary \ref{co.-motZ}, we can express the motive $h^i(Z)$ as
$ \bigoplus_{j+k=i}m_j.\Le^j\otimes h^k(Y)$.
Here $m_j$ is the number of $j$-codimensional cells on $\F$.
Hence
 $$L^{n-i}_\al=\oplus_{j+k=i}m_j. L^{f-j}_\gamma\otimes L_\beta^{m-k}$$
and 
$$L_\al^{n-i}:h^i(Z)\lrar h^{2n-i}(Z)$$
gives an isomorphism.
A similar argument holds for the Lefschetz decomposition.
\end{proof}


\subsection{Varieties with a nef tangent bundle}

In this subsection, we want to deduce the properties discussed in \S 3.2, for varieties with a nef tangent bundle. We want to use the classification results by Campana and Peternell to describe the motive of these vareities.

We first prove the following statement which holds for some wider class of varieties and later deduce some consequences for varieties with nef tangent bundle.

\begin{lemma}\label{cellular}
Suppose $X$ is a nonsingular projective variety such that there is a finite surjective cover
$Z\lrar X$ and $Z$ admits a fibration:
$$
\pi:Z\lrar A
$$
which is a relative cellular variety over an abelian variety $A$. Then $X$ has a 
Chow--K\"unneth decomposition and is finite dimensional.
Furthermore, the motivic Hard Lefschetz Theorem and Lefschetz decomposition holds for $X$.
\end{lemma}

\begin{proof}
Using Lemma \ref{le.-kd} and Corollary \ref{co.-fd}, we conclude that $Z$ has finite dimensional motive and a Chow--K\"unneth decomposition. Now by Lemma \ref{le.-2} and Lemma \ref{le.-3}, we deduce that the variety $X$ has a finite dimensional
motive and a K\"unneth decomposition. Again applying Proposition \ref{pr.-1}, we conclude that $X$ has a Chow--K\"unneth
decomposition.

Since $g$ is a finite morphism, we consider an ample divisor $\al$ on $X'$ which is the pullback of an ample divisor $\eta$ on $X$.
By \cite[Remark 10.6]{Ki}, the motive of
$A$ is finite dimensional and $A$ has a Chow--K\"unneth
decomposition \cite{Sh}. Also the motivic Hard Lefschetz theorem and the Lefschetz decomposition holds for the abelian variety $A$, by Theorem \ref{th.-KuSc}.

Hence, applying Corollary \ref{co.-Hlef},
the Lefschetz operator  $L_\al$ induces isomorphisms 
$$
L^{n-i}_\al: h^i(X')\sta{\simeq}{\lrar} h^{2n-i}(X').
$$
Since the motive $h^i(X)$ is a direct summand of the motive $h^i(X')$, the operator $L^{n-i}_\al$ restricts on 
$h^i(X)$ to the operator $L^{n-i}_\eta$.
This gives a motivic Hard Lefschetz theorem for $X$.
A similar argument also holds for the Lefschetz decomposition.
\end{proof}

\begin{corollary}\label{co.-fCKL}
Suppose $X$ is a nonsingular projective variety with numerically effective tangent bundle such that there is an \'etale cover $X'$ which is  a relative cellular variety over an abelian variety.
 Then the motive of $X$
is finite dimensional and has a Chow--K\"unneth decomposition.
Furthermore, the motivic Hard Lefschetz theorem and Lefschetz decomposition
holds for $X$.
\end{corollary}
\eop

We will see in \S \ref{threefold}, that the assumptions made in the statements of Corollary \ref{co.-fCKL}
are fulfilled by surfaces and threefolds with a nef tangent bundle and hence the above consequences hold for such varieties.

 
\subsection{Chow--K\"unneth projectors }\label{projectors}

Suppose $A$ is an abelian variety with a finite group $G$ acting on $A$.
Then we have

\begin{theorem}\label{th.-ak-jo}
Let $A$ be an abelian variety of dimension $d$ and $G$ be a finite group acting on $A$. Let $h:A\lrar \f{A}{G}$ be the quotient morphism. Suppose $\Delta_A=\sum_{i=0}^{2d}\pi_i$ is the Deninger-Murre Chow--K\"unneth decomposition
for $A$ and let 
$$\eta_i=\f{1}{|G|}(h\times h)_*\pi_i.$$
Then the Chow--K\"unneth decomposition for $\f{A}{G}$ is
$$\Delta_{\f{A}{G}}=\sum_{i=0}^{2d} \eta_i.$$
This Chow--K\"unneth decomposition satisfies Poincar\'e duality, i.e., $\eta_{2d-i}= \,^t\eta_i$, for any $i$.

In addition, $\eta_i$ acts as zero on $CH^j(A/G)_\Q$ for $i<j$ and also for $i>j+d$ in general. In case, $d\leq 4$, we
conclude that $\eta_i$ acts trivially on $CH^j(A/G)_Q$ for $i< j$ and for $i>2j$.
\end{theorem}
\begin{proof}
This is \cite[Theorem 3.1]{De-Mu} when $G$ is trivial and when $G$ is non-trivial this is \cite[Theorem 1.1]{Ak-Jo}.
\end{proof}

Suppose $X$ is a nonsingular projective variety with numerically effective tangent bundle.
By \cite{DPS}, there is a finite \'etale cover 
$$q\,:X'\lrar X$$ such that $X'$ admits a fibration
$$f'\,:X'\lrar A'$$
where $A'$ is an abelian variety and the fibres of $f'$ are Fano varieties with a nef tangent bundle.
We will assume that such Fano varieties are rational homogenous spaces (see \S \ref{remarks}).

Since the fibers $F$
of $f'$ are simply connected and are rationally
connected varieties \cite{Kol}, \cite{KMM}, the homogenous space $F$ does not admit any nonconstant morphisms
to an abelian variety (as there are no nonconstant morphisms from 
${\mathbb P}^1$ to an abelian variety). Hence, 
there is a covering map $q' \, :\, A'\, \longrightarrow\, A$
and a surjective morphism $f\, :\, X\, \longrightarrow\, A$,
such that $q$ is the pullback of $q'$ by $f$.
Now, any morphism between abelian varieties is a composition of a group homomorphism and a translation by an element \cite[Proposition 2.1 a)]{La-Bi}. Since a translation map is an isomorphism and inducing isomorphisms on the Chow groups, without any confusion, we assume that $q'$ is a surjective group homomorphism. Let $G$ be the kernel of $q'$.
Then $A= A'/G$ and $q'\,:A'\lrar A=\f{A'}{G}$.

Moreover, we have a commutative diagram
\begin{eqnarray*}
 X'=A'\times _A X & \sta{q}{\lrar} & X=X'/G \\
\downarrow {f'}\,\,\, & & \downarrow {f} \\
A'\,\,\,\,\,\,\, & \sta{q'}{\lrar} & A 
\end{eqnarray*}
such that $G$ acts on the first factor of $X'$ and has a trivial action on the second factor.

Next, we want to obtain explicit Chow--K\"unneth projectors for the variety $X$ in terms of the projectors of its Albanese reduction.

\begin{proposition}\label{pr.-ckx}
Let $\Delta_A=\sum_i \pi^A_i$ be the Deninger-Murre Chow--K\"unneth decomposition for $A$. Then 
$\Delta_X=\sum_i\pi^X_i$ where $$\pi^X_i=\sum_{j+k=i}\tilde\pi_k.(f\times f)^*\pi^A_j$$ is a Chow--K\"unneth decomposition for
$X$. Here $\{\tilde\pi_k\}$ are the images of the projectors defined on 
the cover $X'$ as in Lemma \ref{co.-motZ}.
 Here $d=\m{dim }A$ and $n=\m{dim }F$.
\end{proposition}
\begin{proof}
Firstly, we claim that the (rational) pushforward of the Deninger-Murre Chow--K\"unneth decomposition $\Delta_{A'}$, is the Deninger-Murre Chow--K\"unneth decomposition on $A$;
by Lemma \ref{le.-2}, the (rational) pushforward of a Chow--K\"unneth decomposition of $\Delta_{A'}$ is a K\"unneth decomposition.
Consider the multiplication map, for any $n\in \Z$, 
$$id \times n :A'\times A' \lrar A'\times A'$$  
There is an eigenspace decomposition of the rational Chow groups of $A'\times A'$. These eigenspaces are defined in \cite[Corollary 2.21]{De-Mu}.
We further apply \cite[Proposition 2 c)]{Be} to obtain that the pushforward map $(q'\times q')_*$ respects the eigenspace decomposition of the Chowgroups of $A' \times A'$ and $A\times A$.
In particular, consider the Deninger--Murre Chow--K\"unneth decomposition
$$
\Delta_{A'}= \sum_{i=0}^{2d}\pi_i
$$
such that $(id \times n)^*\pi_i=n^i\pi_i$, for all $n\in \Z$.
Furthermore, this is a unique decomposition with these properties \cite[Theorem 3.1]{De-Mu}.

Then the (rational) pushforward to $A\times A$ of the above decomposition of $\Delta_{A'}$ also satisfies the same properties and by the uniqueness of such a decomposition gives the Deninger--Murre Chow--K\"unneth decomposition of $A$.

Consider the Chow--K\"unneth decomposition of $X'$, using Corollary 
\ref{co.-motZ} :
$$
(X',\Delta_{X'})= \bigoplus_i(X',\pi^{X'}_i)
$$
where $(X',\pi^{X'}_i)= \bigoplus_{j+k}m_j.\Le^j\otimes (A',\pi^{A'}_k)$.
Here $m_j$ is the number of $j$-codimensional cells on a fibre $\F$ and $\pi^{A'}_k$ are the Deninger--Murre Chow--K\"unneth projectors of $A'$.

Now consider the (rational) pushforward $(q\times q)_*$ of this Chow--K\"unneth decomposition of $X'$ which is now a K\"unneth decomposition of $X$, by Lemma \ref{le.-2}. Since the pushforward of the projectors $\pi^{A'}_k$ are precisely the Deninger--Murre
Chow--K\"unneth projectors of $A$ (from the above discussion), we deduce that the pushforward of the projectors $\pi_i^{X'}$ satisfy the orthogonality properties. 

\end{proof}

Suppose $X'\lrar A'$ is a relative cellular space.

Then there is a sequence of closed embeddings
$$\emptyset=Z_{-1}\subset Z_{0} \subset...\subset Z_n=X'$$
such that $\pi_k:Z_k\lrar Y$ is a flat projective $A'$-scheme. Further, for any $k=0,1,...,n$, the open
complement $Z_k-Z_{k-1}$ is $A'$-isomorphic to an affine space $\A^{m_k}_{A'}$ of relative dimension $m_k$.
Denote $i_k:Z_k\hookrightarrow Z$.

Then we have
\begin{lemma}\label{le.-hcg} 
For any $a,b\in \Z$, the map
$$\bigoplus_{k=0}^nH_{a-2m_k}(A',b-m_k)\lrar H_a(X',b)$$
$$(\al_0,...,\al_n)\mapsto \sum_{k=0}^n(i_k)_*\pi_k^*\al_k$$
is an isomorphism.
Here $H_a(A',b)=CH_b(A',a-2b)$, the higher Chow groups of $A'$.
\end{lemma} 
\begin{proof}
See \cite[Theorem, p.371]{Ko}.
\end{proof}

\begin{remark}\label{re.-chow}
The above Lemma \ref{le.-hcg} can equivalently be restated to express the Chow groups of $X'$ as 
$$CH^r(X')_\Q= \bigoplus_{k=0}^r (\oplus_\al\Q[\omega_k^\al]).f^*CH^k(A')_\Q,$$
Here $\omega_k^\al$ are the $r-k$ codimensional relative cells and $\al$ runs over the indexing set of $r-k$ codimensional relative cells in the $A'$-scheme 
$X'$.
\end{remark}

\begin{lemma}\label{le.-chx}
Suppose $X$ is a nonsingular projective variety and has a finite \'etale cover
$X'\lrar A$ which is a relative cellular variety over an abelian variety $A$. 
Then the Chow groups of $X$ can be written as
$$
CH^r(X)_\Q\,=\, \bigoplus_{k=0}^r (\sum_\al\Q.[\ov\omega_k^\al]).f^*CH^k(A)_\Q.
$$
Here $\ov\omega_k^\al$ are the images of the $r-k$ codimensional relative cells $\omega_k^\al$ and $\al$ runs over the indexing set of $r-k$ codimensional relative cells in the $A'$-scheme $X'$ (see Remark \ref{re.-chow}).
\end{lemma}
\begin{proof}
Using \cite[Example 1.7.6, p.20]{Fu}, the Chow groups of $X$ are written as
\begin{eqnarray*}
CH^r(X)_\Q & = & (CH^r(X')_\Q)^G ,\m{ the group of }G-\m{invariants in }CH^r(X')_\Q \\
           & = & \bigoplus_{k=0}^r (\sum_\al\Q[\omega_k^\al]).({f'}^*CH^k(A')_\Q)^G , \m{ by Remark }\ref{re.-chow}\\
           & = & \bigoplus_{k=0}^r (\sum_\al\Q[\omega_k^\al]).{f'}^*(CH^k(A')_\Q^G) \\
           & = & \bigoplus_{k=0}^r (\sum_\al\Q[\ov\omega_k^\al]).f^*CH^k(A)_\Q ,\m{ again by loc.cit.}
\end{eqnarray*}
\end{proof}

\begin{corollary}\label{Murre}
The Chow--K\"unneth decomposition obtained in Proposition \ref{pr.-ckx} satisfies Poincar\'e duality, i.e., $\pi^X_{2d-i}= \,^t\pi^X_i$, for any $i$.
In addition, $\pi^X_i$ acts as zero on $CH^j(X)_\Q$ for $i<j$ and also for $i>j+d+n$. In case, $d\leq 4$,
we conclude that $\pi^X_i$ acts trivially on $CH^j(X)_Q$ for $i< j$ and for $i>2j$.
\end{corollary}

\begin{proof}
To prove the assertion, notice that the action of the projectors $\pi^X_i$ is given
by the projectors $\tilde{\pi_j}$ (resply. $\pi^A_k$) on $(\sum_\al\Q[\omega_k^\al])$ (resply. on $CH^k(A)_\Q$, for instance see Theorem \ref{th.-ak-jo}).
Since the assertion is true for the projectors $\tilde{\pi_j}$ and $\pi^A_k$ in the desired range, the assertion is true
for the projectors $\pi^X_i$. 
\end{proof}

\subsection{Surfaces and Threefolds with a nef tangent bundle}\label{threefold}
In this subsection we want to show that surfaces and threefolds with a nef tangent bundle satisfy Murre's conjectures.

For this purpose, we recall the classification theorem of such varieties proved by Campana and Peternell.

\begin{theorem}\label{Campana} \cite[Theorem 10.1]{Ca-Pe}
Suppose $X$ is a nonsingular projective threefold. Then the following are equivalent:

1. The tangent bundle $T_X$ is nef.

2. Some \'etale covering $X'$ of $X$ belongs to the following list:

a) $X'=X$ and is a nonsingular Quadric or is the projective space (i.e., $X$ is rational homogenous and $b_2(X)=1$).

b) $X'=X=\p(T_X)$ (i.e., $X$ is rational homogenous and $b_2(X)>1$).

c) $X'=\p(E)$ for a flat rank $3$ vector bundle $E$ on an elliptic curve.

d) $X'=\p(F)\times_C \p(F')$ for flat rank $2$ vector bundles $F$ and $F'$ over an elliptic curve $C$.

e) $X'=\p(E)$ for a flat  rank $2$ vector bundle $E$ on an abelian surface.

f) $X'$ is an abelian threefold.
\end{theorem}

\begin{proposition}\label{threefolds}
Suppose $X$ is a nonsingular projective surface or a threefold which has a nef tangent bundle. Then the motive  of
$X$ is finite dimensional and we have an explicit Chow-K\"unneth decomposition in terms of the projectors of its albanese reduction.
Moreover, Murre's conjectures is fulfilled by these projectors as in Corollary \ref{Murre} and the motivic Hard Lefschetz theorem holds .
\end{proposition} 
\begin{proof}
We only need to show that there is an \'etale cover $X'$ of $X$ which is a relative cellular variety over an
abelian variety.
We now apply Theorem \ref{Campana}  to deduce that there is an \;etale cover $X'$ of $X$ which
is either a projectivization  $\p(E)$ of a vector bundle $E$  over an abelian variety $A$ or a product  of such varieties over an abelian variety. Since the projectivization $\p(E)$ is clearly a relative cellular variety over $A$, the same is true of a product $\p(E)\times_C\p(E')$ over an elliptic curve. A similar classification result holds for surfaces, see \cite[Theorem 3..1]{Ca-Pe}.
Then the proposition  follows from  Corollary \ref{co.-fCKL} and Corollary \ref{Murre}.

\end{proof}

\subsection{Remarks on higher dimensional case}\label{remarks}

Campana and Peternell \cite{Ca-Pe} have made the following conjecture:

\begin{conjecture} \label{conj}
Suppose $X$ is a variety with a nef tangent bundle of dimension $d$. Then there 
is an \'etale cover $X'\lrar X$ such that $X'$ admits a smooth fibration $X'\lrar A$ whose fibres are homogenous spaces and $A$ is an abelian variety.
\end{conjecture}

Furthermore, they also pose the following:

\begin{question} \label{quest}\cite[p.170]{Ca-Pe}:
 Is every Fano manifold with a nef tangent bundle a homogenous variety ? 
\end{question}

The above Conjecture  \ref{conj} is proved in \cite{Ca-Pe} when $d\leq 3$.
In higher dimensions, Demailly-Peternell-Schneider \cite{DPS} prove that the fibres of $X'\rar A$ are Fano manifolds with a nef tangent bundle, if $X$ is as above. An affirmative answer to the above question in all dimensions will prove Conjecture \ref{conj}. If one can show the existence of a relative Chow-K\"unneth decomposition for the fibration
$X'\rar A$ or show that it is a relative cellular variety, then the results of Proposition \ref{threefolds} will hold for
 all varieties which have a nef tangent bundle.


\end{document}